# A fuzzy inventory model considering imperfect quality items with receiving reparative batch and order overlapping

**Hesamoddin Tahami, Hengameh Fakhravar**
Engineering Management & Systems Engineering Department
Old Dominion University
Norfolk, VA 23529, USA
htahami@odu.edu, hfakhrav@odu.edu

**Abstract**

This paper presents an inventory model for imperfect quality items with receiving a reparative batch and order overlapping in a fuzzy environment by employing triangular fuzzy numbers. It is assumed that the imperfect items identified by screening are divided into either scrap or reworkable items. The reworkable items are kept in store until the next items are received. Afterward, the items are returned to the supplier to be reworked. Also, a discount on the purchasing cost is employed as an offer of cooperation from a supplier to a buyer to compensate for all additional holding costs incurred to the buyer. The rework process is error-free. An order overlapping scheme is employed so that the vendor is allowed to use the previous shipment to meet the demand by the inspection period. In the fuzzy model, the graded mean integration method is taken to defuzzify the model and determine its approximation of a profit function and optimal policy. In doing so, numerical examples are rendered to represent the model behavior and, eventually, the sensitivity analysis is presented.

## 1. Introduction

The importance of inventory management systems is growing every day, and many researchers are trying to solve management problems using mathematical models. The economic order quantity (EOQ) model is the basis of advanced inventory systems. By exploring the literature review on inventory systems, it is realized that many efforts have been conducted to provide inventory models in order to eliminate the limitations of the EOQ model. One of the assumptions in the EOQ models is that all the received items are perfect. However, this assumption is not comprehensive for several reasons, including the faulty production process, failure in the process of transportation, etc. So, the effect of imperfect items on inventory systems has become one of the interesting topics for many researchers to provide more practical models. Porteus (1986), followed by Rosenblatt and Lee, presented the significant connection between imperfect quality and lot sizing. Schwaller (1988) assumed that the imperfect items received in the lot would result in the inspection cost. In addition, Zhang and Gerchak (1990) studied an EOQ model with the effect of a joint lot sizing and screening, in which the imperfect items were random variables. Further, Salameh and Jaber (2000) investigated an economic production quantity model for defective items with a known probability distribution. Therefore, they assumed that, by the end of the inspection time, the imperfect items were sold as a single batch. In the same year, Cardenas-Barron (2000) did not deviate from the main idea, but pointed out and rectified an existing error in the model, which was devised by Salameh and Jaber. In (2002), a simple method was proposed by Cardenas-Barron and Goyal to find the optimum value of production quantity for the model given by Salameh and Jaber. Subsequently, Papachristos and Konstantaras (2006) examined the imperfect inventory model, given that the imperfect items were random variables. Jaber et al. (2011) extended the inventory model presented by Salameh and Jaber (2000), while considering the effect of error in the screening procedure.

Zanoni, Jaber, and Zavanella (2013) further presented an EOQ inventory model considering imperfect quality items. In the same year, Moussawi-Haidar, Salameh, and Nasr (2013) suggested an inventory model, in which lot-sizing, defective items, quality control were combined. Karimi-Nasab and Sabri-Laghaie (2014) formulated a new imperfect production inventory model, in which the imperfect items were randomly produced. Currently, Moussawi-Haidar, Salameh, and Nasr (2014) considered the effect of imperfect items and deterioration. It should be noted that, in all the above-mentioned models, no shortage has been assumed during the inspection process, which was based on the model by Salameh and Jaber (2000). Since in several successive inspection processes, defective items may have been found to cause a shortage, the supposed assumption could not be correct. This fault was discussed by Papachristos and Konstantaras (2006), who concluded that the simple formula could not be found to prevent the occurrence of shortage during the inspection process. Luckily, Maddah et al. (2010) developed a pragmatic method to overcome this fault. This method, called "an order overlapping scheme," lets the vendor use the previous order to meet the demand during the inspection process. This new approach can effectively prevent the occurrence of shortage during the inspection process. Therefore, this idea was incorporated into our model.

Another unrealistic assumption considered in the above models was that imperfect goods could just be sold at their salvage value and could not be reworked. However, many researchers have discovered this fault and incorporated the idea of reworking a part of imperfect items into their models. For example, Salameh and Hayek (2001) proposed an economic production quantity model, in which the defective items were reworked by the end of the production time. Yu et al. (2012) studied an EOQ model, in which a part of defective items could be used as good items. Treviño-Garza et al. (2014), Ouyang et al. (2014), and many other researchers have also investigated the effect of the reworked process on the inventory models. It should be noted that the above inventory models consider that the reworkable items are sent back to be reworked and returned as the perfect items through the same period; however, in our model, we assumed that reworkable items were kept in the buyer's warehouse until the next shipment arrived. Then, the supplier replaced the reworkable items with the perfect ones and sent them within the next order before the current lot was used up. In the present paper, it was assumed that the following lot was received from the supplier as "reparative batch". Also, in the previous papers, it was assumed that the perfect item holding costs and scrap item holding costs were the same. However, Wahab and Jaber (2010) and Tahami et al. (2016) presented an imperfect EOQ inventory model with different holding costs and learning in the inspection.

In all the above models, researchers have only considered all the parameters and variables as crisp values. Although crisp models offer an overview of the approach of inventory systems under various assumptions, they are not able to provide factual terms. As a result, exerting crisp models in general, can lead to errors in decision-making. Moreover, in crisp models, inventory managers must be flexible in determining the economic lot size to cause uncertainty cost reduction. Furthermore, the use of fuzzy systems for solving inventory problems, rather than using probability systems, generates more appropriate solutions. Fuzzy sets introduced by Zadeh (1965) drew the attention of many researchers in the inventory management topics.

Sommer (1981) developed a fuzzy scheduling inventory model considering a constraint in warehouse capacity. Stanieski and Kacprzyk (1982) studied a long-term inventory method applying fuzzy decision-making models. Park (1987) presented an EOQ model for interpreting a fuzzy set theory. Chang et al. (1998) developed an EOQ inventory model considering the backorder as a triangular fuzzy number. Yao and Chiang (2003) provided an inventory model without backorder, considering the fuzzy storing cost defuzzified by centroid and signed distances. Chang (2004) proposed an imperfect inventory model considering the fuzzy annual demand and fuzzy imperfect rate. Mahata (2011) studied a fuzzy EOQ inventory model with two-phase trade credits for deteriorating items in the fuzzy sense. Since then, Tahami et al (2019) has made significant contributions to controllable lead-time literature.

Bjork and Carlsson (2005) fuzzified lead time components and studied the effect of flexibility in lead time on the distributors. Bjrk (2008) investigated an EPQ model without any shortages by fuzzifing the decision variable and cycle time. Furthermore, Bjrk (2009) proposed the EOQ model considering triangular fuzzy numbers for the demand and lead time.

As it is obvious from the above-mentioned literature, none of the authors has presented an imperfect EOQ inventory model, either scrap or re-workable, along with receiving reparative batch considering various holding costs for perfect and scrap items under fuzzy conditions in the model parameters. Therefore, we tried to eliminate the gap in the literature. In this paper, scrap items were being sold for salvage value by the end of the inspection period. Upon the completion of the screening process, the buyer notifies the supplier of the number of reworkable items; however, unlike some of the previous articles, here, it is assumed that reworkable items are stored in the buyer's warehouse until the next shipment arrives. Then, the supplier replaces the reworkable items with the perfect ones and sends them within the next order before the current lot is exhausted. Totally, the major distinction between this paper and others lies in fuzziness in the model parameter, the various assumptions on imperfect items, employing an overlapping scheme to prevent shortages during the inspection period, discount rate provision of the purchasing cost to maintain a cooperative relationship, and considering receiving reparative.

The structure of this paper is as follows. The problem statement is given in the second section. Next, the mathematical model is presented. Afterward, numerical examples and sensitivity analysis are given. Finally, the conclusion section is provided.

## 2. Problem Statement

In this section, the problem is introduced with more details. An imperfect EOQ inventory model is presented. All the items received on a shipment are required to be inspected. The imperfect items that are identified through screening are divided into either scrap or reworkable items. By the end of the inspection period, the scrap items are sold at a price of salvage value. Then, the buyer declares the number of reworkable items; however, unlike some of the previous articles in which reworkable items are assumed to be sent back to the supplier and returned as the perfect items through the same period, the proposed model is assumed that reworkable items are kept in a buyer's warehouse until the next shipment arrives. Then, the supplier replaces the reworkable items with the perfect ones and sends them within the next order before the current lot is exhausted. By doing so, the supplier's costs (e.g., transportation costs) are reduced and, instead, the buyer's costs (e.g., holding costs) are raised. As a result, a coordinated policy should be employed so that economic benefits can be provided for both the buyer and the supplier. Discount on purchase costs can be used as an offer of cooperation from supplier to buyer (i.e., the discount compensates for all additional holding cost incurred to the buyer). Moreover, to eliminate shortages within the inspection period, an "overlapping scheme" is employed: similar to Maddah et al.'s (2010) idea that let the buyer supply his/her needs from the previous order during the inspection process. Also, it is assumed that the holding costs for scrap items and perfect items are not the same. Besides, the input parameter $D$ is considered a triangular fuzzy number and applies a graded mean integration method as a defuzzification method to obtain the optimum values. In addition, Yahoodik et al. (2020) and Tahami and Fakhravar (2020) stated that demand is stochastically distributed in its nature in most industries.

Following are the assumptions considered in this paper:
- Item demand is constant over time.
- The input parameter $D$ is the triangular fuzzy number.
- A graded mean integration method is applied as defuzzification so that the optimum value of the profit function in the fuzzy case could be found.
- Shortages are not allowed.
- The holding cost for reworkable items is different from and higher than the holding cost for scrap items.
- A discount on the purchasing cost is applied to make up for the extra holding cost belonging to the buyer.
- An order overlapping scheme is incorporated into the model.
- The demand and screening processes proceed concurrently, but $D < x$.
- Reworking process is error-free.

## 3. Mathematical Modeling

The notations used in the paper are as follows:

$\tilde{D}$ : Demand per year, nonnegative triangular fuzzy number with parameter ($q$, $r$, $s$)

$x$ : Inspection rate

$A$ : Ordering cost per cycle

$r_s$ : Percentage rate of scrap items (random variable)

$r_w$ : Percentage rate of reworkable items (random variable)

$f(r_s)$: $r_s$ Probability density function

$f(r_w)$: $r_w$ Probability density function

$s$ : Selling value per unit

$w$ : Salvage value per unit

$d$ : Unit inspection cost

$h_w$ : Reworkable or perfect item holding cost rates per unit per cycle

$h_s$ : Scrap item holding cost rate per unit per cycle

$\alpha$ : Discount rate for procurement cost

$c$ : Purchasing cost per unit

$t_1$ : Screening length per cycle

$T$ : Length of cycle

$H_s(Q)$ : Scrap item holding cost per cycle

$H_w(Q)$ : Perfect or reworkable item holding costs per cycle

$TP(Q)$ : Total profit per cycle

$TPU(Q)$ : Net profit per unit time

$Q$ : Order size per cycle (decision variable)

Considering that the demand rate $D$ is a fuzzy number; however, other components of the model are all crisp constant. We represent the demand rate by a triangular fuzzy number as given below:

$$\tilde{D} = (q, r, s) \qquad (1)$$

And the membership function is as follows:

$$\mu_{\tilde{D}}(x) = \begin{cases} \dfrac{x-q}{r-q} & \text{if } q \leq x \leq r \\ \dfrac{s-x}{x-r} & \text{if } r \leq x \leq s \\ 0 & \text{otherwise} \end{cases} \qquad (2)$$

Fig. 1 presents the behavior of the proposed model per cycle. The 100% inspection process is finished at time $t_1$. To avoid shortages, the overlapping scheme is used and it is supposed that the demand by the screening time is at least the same as the number of perfect quality items. It means that, for $0 \leq t \leq t_1$:

$$xt_1(1 - r_s - r_w) \geq \tilde{D} t_1 \qquad (3)$$

which yields:

$$x \geq \dfrac{\tilde{D}}{(1 - r_s - r_w)} \qquad (4)$$

The goal is to obtain $Q$ that maximizes the total profit per year, $TP(Q)$, expressed by:
$$TP(Q) = TR(Q) - TC(Q) \tag{5}$$
Where $TR(Q)$ denotes the revenue per cycle and $TC(Q)$ denotes the total cost per cycle. $TR(Q)$ is obtained through the sale of good items and scrap items, i.e.:
$$TR(Q) = sQ(1-r_s) + \omega Q r_s \tag{6}$$
$TC(Q)$ includes the following four costs:
$$TC(Q) = OC + SC + PC + HC \tag{7}$$
Where $OC$ denotes the ordering cost per cycle ($OC = A$), $SC$ denotes the screening cost per cycle ($SC = dQ$), $PC$ denotes the purchasing cost per cycle ($PC = cQ(1-\alpha)$), and $HC$ denotes the holding cost per cycle, which includes the scrap item holding cost per cycle, $H_s(Q)$, and reworkable or perfect item holding cost per cycle, $H_w(Q)$.

$H_s(Q)$ can be obviously calculated using Fig. 1 as shown in the shaded area:
$$H_s(Q) = h_s \left( \frac{Q^2 r_s}{2x} \right) \tag{8}$$

To compute $H_w(Q)$, the total inventory quantity per cycle should be calculated. According to Fig. 1, it is clear that the sum of the areas of ΔZBC, ΔBGR, GIJR, and ΔRJF minus ΔDEF can express the total inventory quantity per cycle. The area of ΔZBC is the same as that of ΔDEF; therefore, we have:
$$\tilde{V} = \frac{Q^2 r_s}{2x} + \underbrace{\frac{Q^2(1-r_s)}{x}}_{\Box BGR} + \underbrace{\frac{(Q(1-r_s))^2}{2\tilde{D}}}_{\Box RJF} \tag{9}$$

Hence, the holding cost $\tilde{H}_w(Q)$ is as follows:
$$\tilde{H}_w(Q) = h_w \times \tilde{V} = h_w \left( \frac{Q^2 r_s}{2x} + \frac{Q^2(1-r_s)}{x} + \frac{(Q(1-r_s))^2}{2\tilde{D}} \right) \tag{10}$$

Thus:
$$T\tilde{C}(Q) = A + dQ + cQ\left(1 - \frac{r_w Q}{\tilde{D}}\right) + h_s\left(\frac{Q^2 r_s}{2x}\right) + h_w\left(\frac{Q^2 r_s}{2x} + \frac{Q^2(1-r_s)}{x} + \frac{(Q(1-r_s))^2}{2\tilde{D}}\right) \tag{11}$$

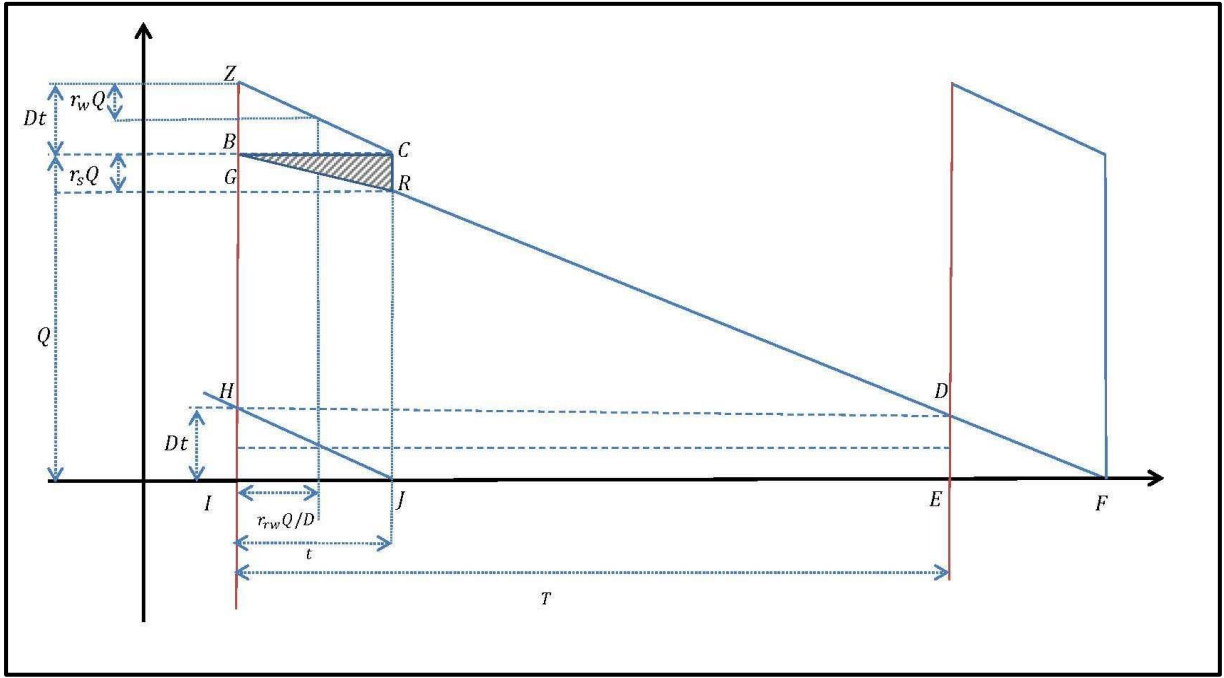

Fig. 1. Inventory model

Through items simplification, the expression for total cost per cycle can be calculated by:

$$T\tilde{C}(Q) = A + dQ + cQ\left(1 - \frac{r_w Q}{\tilde{D}}\right) + (h_s + h_w)\left(\frac{Q^2 r_s}{2x}\right) + h_w\left(\frac{Q^2(1-r_s)}{x}\right) + h_w\left(\frac{(Q(1-r_s))^2}{2\tilde{D}}\right) \quad (12)$$

By substituting Eqs. (12) and (6) in Eq. (5), the total profit per cycle is obtained by:

$$T\tilde{P}(Q) = sQ(1-r_s) + \omega Q r_s - A - dQ - cQ\left(1 - \frac{r_w Q}{\tilde{D}}\right) - (h_s + h_w)\left(\frac{Q^2 r_s}{2x}\right) \\ - h_w\left(\frac{Q^2(1-r_s)}{x}\right) - h_w\left(\frac{(Q(1-r_s))^2}{2\tilde{D}}\right) \quad (13)$$

Furthermore, it is considered that the expected value of $T\tilde{P}(Q)$ (i.e., $E[T\tilde{P}(Q)]$) is calculated, in which the expected values ($E[1-r_s]$, $E[r_s]$, and $E[r_w]$) are used instead of $1-r_s$, $r_s$, and $r_w$, respectively. The expected net profit per unit time is calculated by applying the renewal reward theorem (Ross, 1996) (i.e., dividing $T\tilde{P}(Q)$ by the cycle length, $\tilde{T} = \frac{(1-r_s)Q}{\tilde{D}}$) as follows:

$$E[T\tilde{P}U(Q)] = \frac{\tilde{D}\left(s(1-E(r_s)) + \omega E(r_s) - c\left(\frac{1-E(r_w)Q}{\tilde{D}}\right) - d\right) - \left(\frac{A\tilde{D}}{Q}\right)}{1 - E(r_s)} \\ - \frac{Q}{2(1-E(r_s))}\left(\frac{\tilde{D}(2h_w - h_w E(r_s) + h_s E(r_s))}{x} + h_w\left(E(1-r_s)^2\right)\right) \quad (14)$$

Now, we consider:

$$u = \frac{1}{1-E(r_s)}\left\{s(1-E(r_s))+\omega E(r_s)-c-d-\frac{A}{Q}\right\}$$

$$W = \frac{Q}{2(1-E(r_s))}\left\{\frac{(2h_w - h_w E(r_s) + h_s E(r_s))}{x}\right\} \tag{15}$$

$$X = \frac{h_w\left(E(1-r_s)^2\right)}{2(1-E(r_s))}$$

Hence, by substituting Eq. (15) in Eq. (14), the annual fuzzy net profit function is illustrated by:

$$E(\tilde{T}PU(Q)) = \tilde{D}(u) + \frac{cE(r_w)}{1-E(r_s)}Q - \tilde{D}(W) - Q(X) \tag{16}$$

Therefore, the annual fuzzy net profit function is illustrated by a nonnegative triangular fuzzy number as follows:

$$E(\tilde{T}PU(Q)) = (a_1, b_1, c_1) \tag{17}$$

Where $a_1, b_1, c_1$ can be obtained below. According to Eq. (16), we have:

$$E(\tilde{T}PU(Q)) - \frac{cE(r_w)}{1-E(r_s)}Q + Q(X) = \tilde{D}(u - W)$$

$$\Rightarrow \tilde{D} = \frac{E(\tilde{T}PU(Q)) - \left(\frac{cE(r_w)}{1-E(r_s)}\right)Q + Q(X)}{u - W} \tag{18}$$

Using Eqs. (18) and (17), we have:

$$\text{for } \tilde{D} \le a \Rightarrow \frac{E(\tilde{T}PU(Q)) - \left(\frac{cE(r_w)}{1-E(r_s)}\right)Q + Q(X)}{u - W} \le a$$

$$\Rightarrow E(\tilde{T}PU(Q)) \le \underbrace{a(u - W) + \frac{cE(r_w)}{1-E(r_s)}Q - Q(X)}_{a_1} \tag{19}$$

Also:

$$\text{for } a \le \tilde{D} \le b \Rightarrow a \le \frac{E(\tilde{T}PU(Q)) - \left(\frac{cE(r_w)}{1-E(r_s)}\right)Q + Q(X)}{u - W} \le b$$

$$\Rightarrow \underbrace{a(u-W) + \frac{cE(r_w)}{1-E(r_s)}Q - Q(X)}_{a_1} \le E(\tilde{T}PU(Q)) \le \underbrace{b(u-W) + \frac{cE(r_w)}{1-E(r_s)}Q - Q(X)}_{b_1} \tag{20}$$

and

for $b \leq \tilde{D} \leq c \Rightarrow b \leq \dfrac{E(\tilde{TPU}(Q)) - \left(\dfrac{cE(r_w)}{1-E(r_s)}\right)Q + Q(X)}{u - W} \leq c$

(21)

$\Rightarrow \underbrace{b(u-W) + \dfrac{cE(r_w)}{1-E(r_s)}Q - Q(X)}_{b_1} \leq E(\tilde{TPU}(Q)) \leq \underbrace{c(u-W) + \dfrac{cE(r_w)}{1-E(r_s)}Q - Q(X)}_{c_1}$

Thus, $E(\tilde{TPU}(Q))$ is a triangular fuzzy number with three points $(a_1, b_1, c_1)$ as follows:

$$E(\tilde{TPU}(Q)) = \begin{cases} a_1 = a(u-W) + \dfrac{cE(r_w)}{1-E(r_s)}Q - Q(X) \\ b_1 = b(u-W) + \dfrac{cE(r_w)}{1-E(r_s)}Q - Q(X) \\ c_1 = c(u-W) + \dfrac{cE(r_w)}{1-E(r_s)}Q - Q(X) \end{cases}$$

(22)

$E(\tilde{TPU}(Q))$ is deffuzified by employing the graded mean integration method using Formula (A. 1) (see Appendix) as follows:

$$P(E(\tilde{TPU}(Q))) = \dfrac{1}{6}(a_1 + 4b_1 + c_1) = \dfrac{1}{6}\begin{pmatrix} \left\{a(u-W) + \dfrac{cE(r_w)}{1-E(r_s)}Q - Q(X)\right\} + \\ 4\left\{b(u-W) + \dfrac{cE(r_w)}{1-E(r_s)}Q - Q(X)\right\} + \\ \left\{c(u-W) + \dfrac{cE(r_w)}{1-E(r_s)}Q - Q(X)\right\} \end{pmatrix}$$

(23)

The target is to maximize $P(E(\tilde{TPU}(Q)))$. Because $P(E(\tilde{TPU}(Q)))$ is concave at $Q$, $\dfrac{\partial^2 E(\tilde{TPU}(Q))}{\partial^2 Q} = \dfrac{-2A\tilde{D}}{Q^3(1-E(r_s))} < 0$, then the optimum lot size $Q^*$ can be calculated by differentiating $P(E(\tilde{TPU}(Q)))$ with respect to $Q$ and setting the partial derivatives is equal to zero.

$$\dfrac{\partial P(E(\tilde{TPU}(Q)))}{\partial Q} = \dfrac{1}{6}(a + 4b + c)\left(\dfrac{\partial u}{\partial Q} - \dfrac{\partial W}{\partial Q}\right) + \dfrac{cE(r_w)}{1-E(r_s)} - X = 0 \quad (24)$$

where

$$\frac{\partial u}{\partial Q} = \frac{\partial}{\partial Q}\left(\frac{\left\{s(1-E(r_s))+\omega E(r_s)-c-d-\frac{A}{Q}\right\}}{1-E(r_s)}\right) = \frac{A}{Q^2} \times \frac{1}{1-E(r_s)}$$

$$\frac{\partial W}{\partial Q} = \frac{\partial}{\partial Q}\left(\frac{Q\left\{\frac{(2h_w - h_w E(r_s) + h_s E(r_s))}{x}\right\}}{2(1-E(r_s))}\right) = \frac{\frac{(2h_w - h_w E(r_s) + h_s E(r_s))}{x}}{2(1-E(r_s))} \tag{25}$$

Hence, substituting Eq. (25) in Eq. (24):

$$\frac{1}{6}(a+4b+c)+(\frac{A}{Q^2}\frac{1}{1-E(r_s)}-\frac{2h_w - h_w E(r_s)+h_s E(r_s)}{2x(1-E(r_s))})+\frac{cE(r_w)}{1-E(r_s)}-\frac{h_w\left(E(1-r_s)^2\right)}{2(1-E(r_s))}=0 \tag{26}$$

By simplifying Eq. (26), Eq. (27) can be obtained by:

$$\frac{A}{Q^2} = \frac{-6cE(r_w)}{(a+4b+c)} + \frac{3h_w\left(E(1-r_s)^2\right)}{(a+4b+c)} + \frac{\left\{\frac{(2h_w-h_w E(r_s)+h_s E(r_s))}{x}\right\}(a+4b+c)}{2(a+4b+c)} \tag{27}$$

This yields to:

$$Q^* = \sqrt{\frac{2A(a+4b+c)}{(a+4b+c)\left\{\frac{(2h_w-h_w E(r_s)+h_s E(r_s))}{x}\right\}-12cE(r_w)+6h_w\left(E(1-r_s)^2\right)}} \tag{28}$$

The optimum annual total profit $P\left(E\left(\tilde{T}PU(Q)\right)\right)$ is obtained by the direct substitution of Eq. (28) in Eq. (23). Note that, when the input parameter $D$ is a real number, that is $a=b=c=D$, when screening rate is large enough, that is the inspection process is finished simultaneously by the receiving an order, and finally when items are categorized as only perfect or imperfect (no reworkable items so that no discount on purchasing cost), $Q^*$ in Eq. (28) is equivalent to:

$$Q^* = \sqrt{\frac{2AD}{h_w\left(E(1-r_s)^2\right)}} \tag{29}$$

which is the same as the results by Shih (1980) and Silver (1976). It shows that the proposed model is accurate. In addition, it should be noted that, when the input parameter $D$ is a real number, that is $a=b=c=D$, and if all items are assumed to be perfect, our model becomes an equivalent to the EOQ inventory model.

## 4. Numerical Study

In this section, the behavior of our model is investigated by applying numerical examples, and the impact of applying fuzzy case into the model is also investigated. Assume the following values and the input parameters for an inventory model in the crisp case:

$A = 100$ $ per cycle
$D = 50,000$ per unit per year
$x = 175200$ per year per unit
$h_w = 5$$ per year per unit
$s = 50$$ per unit
$d = 0.5$$ per unit
$w = 20$ per unit
$h_s = 2$$ per year per unit
$c = 25$$ per unit

Also, $E(r_s) = 0.02 \quad E(r_w) = 0.05 \quad E[(1-r_s)^2] = 0.9605$

The optimum lot size $Q^*$ and the optimal annual total profit $E[TPU(Q)]$ of a crisp case, in which $a=b=c=D=50000$, can be derived easily from Eqs. (28) and (23), respectively. It is obtained that:

$Q^* = 1395 \quad E[TPU(Q)]^* = 1212072$

Some triangular fuzzy numbers are assigned for the input parameter $\tilde{D}$ in Table 1 to illustrate the fuzzy model developed in Section 3. Then, by using the GMI method, the defuzzified values are specified. The defuzzified values and the corresponding percentage difference from the crisp values (denoted by $p_D$ for the component $\tilde{D}$) are also shown in Table 1. For each set of triangular fuzzy numbers, the optimal lot size $Q^*$ and the optimal annual total profit $P(E(\tilde{TPU}(Q)))$ are derived from Eqs. (28) and (23). The findings are summarized in Table 2. This table represents variations in $Q^*$ and annual net profit $P(E(\tilde{TPU}(Q)))$ due to fuzziness in parameter $D$. It shows that the optimal values for the expected net profit are fully sensitive to increasing percentage changes in parameter $D$'s fuzziness level, while optimal order quantities are comparatively insensitive to increasing percentage changes in parameter $D$'s fuzziness. Note that the percentage changes in the expected net profit are almost the same as the percentage changes of parameter $D$ at different levels, whereas the lot size order quantity changes marginally. Besides, a one-way sensitivity analysis is conducted to determine the impact of other problem parameters on $Q^*$ and $P(E(\tilde{TPU}(Q)))$. In the numerical examples, the value of one parameter is changed at a time, while the values of the others are not changed. Table 3 shows the values used in the sensitivity analysis for different problem parameters. Then, the optimal order quantity and the annual total net profit are calculated using the values shown in Table 3. The corresponding values are presented in Table 4.



Table 1. Triangular fuzzy numbers for the parameter $\tilde{D}$

| $\tilde{D}$ | $p(\tilde{D})$ | $p_D$ |
|---|---|---|
| (5,000 ; 34,250 ; 68,000) | 35,000 | -30 |
| (12,000 ; 37,500 ; 78,000) | 40,000 | -20 |
| (20,000 ; 45,000 ; 70,000) | 45,000 | -10 |
| (29,000 ; 52,000 ; 93,000) | 55,000 | 10 |
| (42,000 ; 61,000 ; 94,000) | 60,000 | 20 |
| (33,000 ; 61,500 ; 111,000) | 65,000 | 30 |

Table 2. Percentage change in optimum values from the crisp case

| $Q^*$ | % change in $Q^*$ | $P\left(E\left(\tilde{T}PU(Q)\right)\right)$ | % change in $P\left(E\left(\tilde{T}PU(Q)\right)\right)$ |
|---|---|---|---|
| 1277.64 | -8 | 848731.233 | -30 |
| 1322.81 | -5 | 970116.010 | -20 |
| 1361.45 | -2 | 1091503.127 | -10 |
| 1424.23 | 2 | 1334281.969 | 10 |
| 1465.76 | 5 | 1536600.692 | 27 |
| 1473.15 | 6 | 1577064.666 | 30 |

Table 3. Experimental values for the example parameters

| Parameter | Base value | Experimental values | | |
|---|---|---|---|---|
| $x$ | 175200 | 87600 | 175200 | 262800 |
| $h_w$ | 5 | 2.5 | 5 | 10 |
| $h_s$ | 2 | 1 | 2 | 4 |
| $A$ | 100 | 50 | 100 | 200 |
| $d$ | 0.5 | 0.25 | 0.5 | 1 |
| $s$ | 50 | 25 | 50 | 100 |
| $c$ | 25 | 12.5 | 25 | 50 |
| $w$ | 20 | 10 | 20 | 40 |
| $E(r_s)$ | 0.02 | 0.01 | 0.02 | 0.03 |
| $E(r_w)$ | 0.05 | 0.025 | 0.05 | 0.075 |

Table 4. Order quantity and expected net profit per unit time for the experimental values

| $x$ | $Q$ | $P\left(E\left(\tilde{T}PU(Q)\right)\right)$ | $h_w$ | $Q$ | $P\left(E\left(\tilde{T}PU(Q)\right)\right)$ |
|---|---|---|---|---|---|
| 87600 | 1119.7 | 1210274.6 | 2.5 | 2746.9 | 1215672.9 |
| 175200 | 1394.9 | 1212072 | 5 | 1394.9 | 1212072 |
| 262800 | 1544.2 | 1212779.7 | 10 | 885 | 1207858 |
| $A$ | $Q$ | $P\left(E\left(\tilde{T}PU(Q)\right)\right)$ | $d$ | $Q$ | $P\left(E\left(\tilde{T}PU(Q)\right)\right)$ |
| 50 | 986.4 | 1214215.1 | 0.25 | 1394.9 | 1224827.6 |
| 100 | 1394.9 | 1212072 | 0.5 | 1394.9 | 1212072 |



| 200 | 1972 | 1209042.5 | 1 | 1394.9 | 1186562.3 |
| --- | --- | --- | --- | --- | --- |
| $h_s$ | $Q$ | $P\left(E\left(\tilde{T}PU(Q)\right)\right)$ | $s$ | $Q$ | $P\left(E\left(\tilde{T}PU(Q)\right)\right)$ |
| 1 | 1395.7 | 212076.6 | 25 | 1394.9 | -37927.4 |
| 2 | 1394.9 | 1212072 | 50 | 1394.9 | 1212072 |
| 4 | 1393.4 | 1212064.4 | 100 | 1394.9 | 3712072.5 |
| $c$ | $Q$ | $P\left(E\left(\tilde{T}PU(Q)\right)\right)$ | $w$ | $Q$ | $P\left(E\left(\tilde{T}PU(Q)\right)\right)$ |
| 12.5 | 1251.1 | 1848986.5 | 10 | 1394.9 | 1201868.5 |
| 25 | 1394.9 | 1212072 | 20 | 1394.9 | 1212072 |
| 50 | 1946.5 | -61364.6 | 40 | 1394.9 | 1232480.7 |
| $E(r_w)$ | $Q$ | $P\left(E\left(\tilde{T}PU(Q)\right)\right)$ | $E(r_s)$ | $Q$ | $P\left(E\left(\tilde{T}PU(Q)\right)\right)$ |
| 0.025 | 1251.1 | 121123.3 | 0.01 | 1393.8 | 1214974.9 |
| 0.05 | 1394.9 | 1212072 | 0.02 | 1394.9 | 1212072 |
| 0.075 | 1603.5 | 1213024 | 0.03 | 1396.1 | 1209110.4 |

Figs. 2 and 3 display a tornado diagram as a graphical result of the sensitivity analysis. These represents how the order quantity and the annual total net profit are changing, while the model parameters are independently varying from their low value to the high ones. The length of each bar in the diagram shows the extent to which the optimal order quantity and the annual net profit are sensitive to the bar's corresponding model parameter. It can be observed from Fig. 2 that the model's parameters with the greatest impact on the optimum order size is $h_w$. As other parameters have their base values, the perfect or reworkable item holding cost rate per cycle is differing from 2.5 to 10, while the value of the order quantity changes from 2746.9 to 885. This finding shows that a larger amount of $h_w$ can highly affect the order quantity. Moreover, it can be observed from Fig. 3 that the model's parameters with the greatest impact on the annual net profit are the unit screening cost. As other parameters have their base values, when $d$ varied from 0.25 to 1, the value of the annual net profit changes from 1224827 to 1186562. As a result, the values of these parameters should be carefully estimated, because they have the most significant impact on the model's cost.

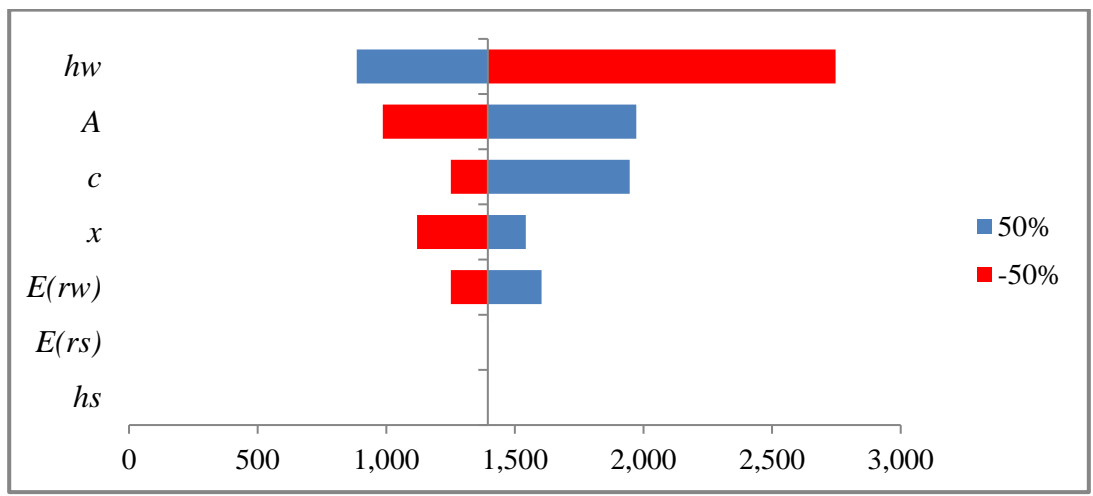

Fig. 2. Tornado diagram for the order quantity



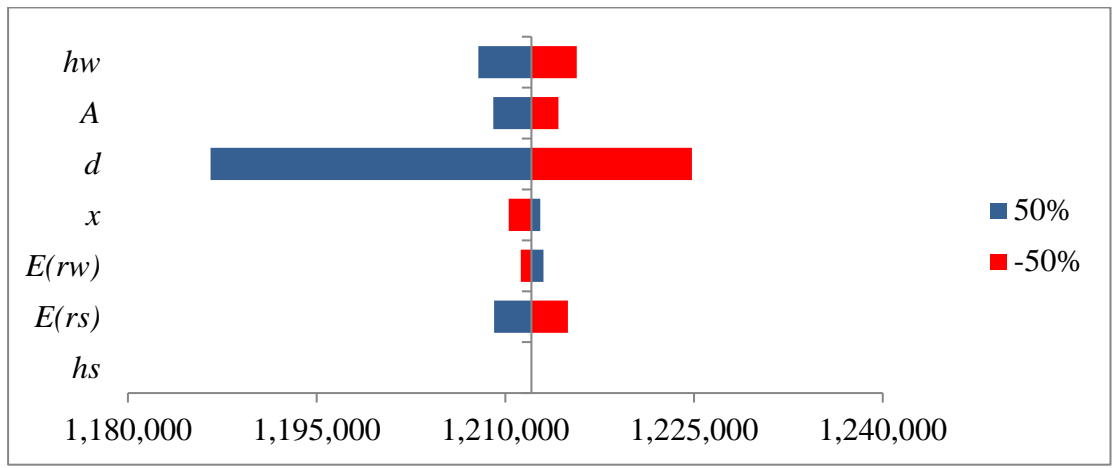

Fig. 3. Tornado diagram for the annual total profit

## 5. Conclusion

As it is known, the input parameters of the EOQ inventory problem cannot be described precisely in a real situation or it may be uncertain because of some uncontrolled factors. Therefore, approximate solution approaches have been represented for the explanation of a series of practical inventory problems. Fuzzy methodologies provide a helpful approach to model ambiguity in human recognition and decision-making. Uncertainties defined by imprecise factors can be illustrated by fuzzy sets. Thus, in the current paper, our goal is to propose the fuzzy inventory model with defective items considering reparative batch and order overlapping. In this model, input parameter ($D$) is considered the fuzzy number to defuzzify the proposed model and determine the approximation of annual profit in the fuzzy sense, we apply the graded mean integration method. Then, the optimal order quantity is calculated to maximize the total profit. The model is solved for triangular fuzzy numbers. It is shown that the EOQ model as well as the models by Silver (1976) and Shih (1980) are just some special cases of our model. In so doing, numerical examples are rendered to represent the model behavior, and then the results of the crisp and fuzzy models are compared with each other. It should be noted that the optimal values of the annual net profit are quite sensitive to increasing percentage changes in the parameter $D$'s fuzziness level, while an optimal order quantity is comparatively insensitive to increasing percentage changes in the parameter $D$'s fuzziness level. The percentage change in the annual net profit is almost the same as the percentage change in the fuzziness level, while the order size changes slightly. A one-way sensitivity analysis is presented to assess the effect of other problem parameters on the order quantity and annual total net profit and to display graphically the sensitivity analysis results as a tornado diagram. To increase the scope of our analysis, the model presented in this paper can be extended in several ways. For example, it can be incorporated with deteriorating items.

**Appendix: Preliminary Concepts in Fuzzy Sets**

**Definition 1** (fuzzy number). A fuzzy set $\tilde{a}$ in the universe of discourse R (set of real numbers) is called a fuzzy number if the following conditions are satisfied:

(1) $\tilde{a}$ is convex.

(2) Membership function $\mu_{\tilde{a}}(x)$, $x \in R$ is at least piecewise continuous.

(iii) $\tilde{a}$ is normal, that is, at least one $x \in R$ exists such that $\mu_{\tilde{a}}(x) = 1$.

**Definition 2** (triangular fuzzy number; TFN). It is a fuzzy number represented with three points as follows:

$$\tilde{v} = (v_1, v_2, v_3)$$

This representation is interpreted as membership function $\mu_{\tilde{v}}(y)$, where:



$$\mu_{\tilde{v}}(y) = \begin{cases} 0 & \text{for } -\infty < y < v_1 \\ \dfrac{y - v_1}{v_2 - v_1} & \text{for } v_1 \leq y < v_2 \\ \dfrac{v_3 - y}{v_3 - v_2} & \text{for } v_2 \leq y \leq v_3 \\ 0 & \text{for } v_3 < y < \infty \end{cases}$$

**Definition 3** (function principle). To perform fuzzy arithmetical operations by TFN, the function principle proposed by Chen (1985) is used. This principle is a suitable method for performing the operations of complex models to prevent arriving at a degenerated solution. This method will be so helpful in handling the fuzzy operations, especially when the crisp model comprises terms of multiple operations of fuzzy numbers. Furthermore, the type of fuzzy membership function will be kept constant during the operations, which helps avoid facing further complexity by arithmetical operations. Now, assume $\tilde{A} = (q_1, r_1, s_1)$ and $\tilde{B} = (q_2, r_2, s_2)$ are two positive TFNs and $\alpha$ are a real number. Employing the functional principle, the operations of the fuzzy numbers $\tilde{A}$ and $\tilde{B}$ are as follows:

i. If $\alpha \geq 0$ $\alpha\tilde{A} = (\alpha q_1, \alpha r_1, \alpha s_1)$ and if $\alpha < 0$ $\alpha\tilde{A} = (\alpha s_1, \alpha r_1, \alpha q_1)$

ii. $\tilde{A} + \tilde{B} = (q_1 + q_2, r_1 + r_2, s_1 + s_2)$

iii. $\tilde{A} - \tilde{B} = (q_1 - s_2, r_1 - r_2, s_1 - q_2)$

iv. $\tilde{A} \times \tilde{B} = (q_1 q_2, r_1 r_2, s_1 s_2)$

v. $\dfrac{\tilde{A}}{\tilde{B}} = (\dfrac{q_1}{q_2}, \dfrac{r_1}{r_2}, \dfrac{s_1}{s_2})$

**Definition 4** (defuzzification). A graded mean integration representation (GMIR) method is used in our paper to transform the fuzzy total profit function to its corresponding crisp function. The main reason is due to the non-linear nature of the function used in this paper. In fact, the fuzzy total profit function in this paper consists of a couple of fuzzy multiplication and division terms; since the GMIR method keeps the shape of the membership function, it is a proper choice to defuzzify the profit function of the model. In the following section, the GMIR method introduced by Chen and Hsieh [4] is described. Assume that $\tilde{A}$ is fuzzy number, $\theta^{-1}$ and $\pi^{-1}$ are the inverse functions of $\theta$ and $\pi$, respectively. The graded mean $\rho$-level value of $\tilde{A}$ is $\dfrac{\rho(\theta^{-1}(\rho) + \pi^{-1}(\rho))}{2}$ and the GMIR of the fuzzy number $\tilde{A}$ is calculated as:

$$\tau(\tilde{A}) = \dfrac{\int_0^{w_A} \dfrac{\rho(\theta^{-1}(\rho) + \pi^{-1}(\rho))}{2} d\rho}{\int_0^1 \rho d\rho} = \int_0^{w_A} (\theta^{-1}(\rho) + \pi^{-1}(\rho)) d\rho$$

where $0 < \rho \leq w_A$ and $0 < w_A \leq 1$.



Now assume $\tilde{A} = (l, m, n)$ is a TFN. The graded mean integration representation of $\tilde{A}$ can be calculated by Formula (A.1), which is as follows:

$$\tau(\tilde{A}) = \frac{\int_0^1 \frac{\rho(\theta^{-1}(\rho) + \pi^{-1}(\rho))}{2} d\rho}{\int_0^1 \rho \, d\rho} = \frac{(l + 4m + n)}{6} \tag{A.1}$$